# Марковские модели макросистем

*Гасников А.В.*
*gasnikov@yandex.ru*

**Задача 1 (модель П. и Т. Эренфестов в форме Кельберта–Сухова)\*.** Рядом стоят две собаки с номерами *1* и *2*. На собаках как-то расположились $N \gg 1$ блох. Каждая блоха в промежутке времени $[t, t+h]$ с вероятностью $\lambda h + o(h)$ ($\lambda = 1$) независимо от остальных перескакивает на соседнюю собаку. Пусть в начальный момент все блохи собрались на собаке с номером *1*. Покажите, что для всех $t \geq \chi \ln N$

$$P\left(\frac{|n_1(t) - n_2(t)|}{N} \leq \frac{3}{\sqrt{N}}\right) \geq 0.99,$$

где $n_1(t)$ – число блох на первой собаке в момент времени $t$, а $n_2(t)$ – на второй (случайные величины). Т.е. относительная разность числа блох на собаках будет иметь порядок малости $O(1/\sqrt{N})$ на больших временах ($T \geq \chi \ln N$). Покажите также, что математическое ожидание времени первого возвращения в начальное состояние (когда все блохи собраны на собаке с номером *1*) равно $2^N/N$.

**Указание.** См. *Narita K.* Asymptotic behavior o continuous-time Markov chain of the Ehrenfest type with applications // Mathematical Models and Stochastic Processes Arising in Natural Phenomena and Their Applications. – 2001.

Введём вектор $p(t) = (p_0(t), ..., p_N(t))^T$, где $p_i(t)$ – вероятность того, что на собаке с номером *1* в момент времени $t \geq 0$ ровно $i$ блох. Тогда в уравнении Колмогорова–Феллера

$$\frac{dp(t)}{dt} = \Lambda^T p(t)$$

инфинитезимальная матрица $\Lambda$ имеет следующий вид:

$$[\Lambda]_{ij} = \begin{cases} 0, & |i - j| > 1 \\ \lambda i, & j = i - 1 \\ \lambda(N - i), & j = i + 1 \\ -\lambda N, & j = i \end{cases}.$$

Описанная марковская динамика имеет закон сохранения числа блох: $n_1(t) + n_2(t) \equiv N$, и это будет единственный закон сохранения. Стационарная (инвариантная) мера имеет вид (теорема Санова):

$$\nu(n_1, n_2) = \nu(c_1 N, c_2 N) = N! \frac{(1/2)^{n_1}}{n_1!} \frac{(1/2)^{n_2}}{n_2!} = C_N^{n_1} 2^{-N} \simeq \frac{2^{-N}}{\sqrt{2\pi c_1 c_2}} \exp(-N \cdot H(c_1, c_2)),$$

где $H(c_1, c_2) = \sum_{i=1}^{2} c_i \ln c_i$. Кстати, сказать, из такого вида стационарной меры, можно получить, что если в начальный момент все блохи находились на одной собаке, то математическое ожидание времени первого возвращения макросистемы в такое состояние будет порядка $2^N / N$. Равновесие данной макросистемы естественно определять как состояние, в малой окрестности которого сконцентрирована стационарная мера (принцип максимума энтропии Больцмана–Джейнса)

$$c^* = \begin{pmatrix} 1/2 \\ 1/2 \end{pmatrix} = \arg\min_{\substack{c_1 + c_2 = 1 \\ c \geq 0}} H(c).$$

Этот же результат можно было получить и при другом порядке предельных переходов (обратном к рассмотренному выше порядку: $t \to \infty$, $N \to \infty$). А именно, сначала считаем, что при $t = 0$ существует предел $c_i(t) \overset{\text{п.н.}}{=} \lim_{N \to \infty} n_i(t)/N$. Тогда (теорема Т. Куртца) этот предел существует при любом $t > 0$, причем $c_1(t)$, $c_2(t)$ – детерминированные (не случайные) функции, удовлетворяющие СОДУ

$$\frac{dc_1}{dt} = \lambda(c_2 - c_1),$$

$$\frac{dc_2}{dt} = \lambda(c_1 - c_2).$$

Глобально устойчивым положением равновесия этой СОДУ будет $c^*$, а $H(c)$ – функция Ляпунова этой СОДУ (убывает на траекториях СОДУ, и имеет минимум в точке $c^*$). Все это можно перенести на общие модели стохастической химической кинетики с (обобщенным) условием детального баланса (см. замечание к задаче 19).

Оценка скорости сходимости макросистемы к равновесию сводится к оценке mixing time соответствующей марковской динамики. Это отдельная задача для данного примера (и моделей, рассмотренных в замечании к задаче 19) может быть решена по-разному. Классическим способом является использование изопериметрического неравенства Чигера, см., например, *Levin D.A., Peres Y., Wilmer E.L.* Markov chain and mixing times. AMS, 2009. Другой способ, еще более ярко подчеркивающий глубокую связь концентрации стационарной меры марковской динамики и оценку mixing time данной динамики, базирующийся на понятии дискретной кривизны Риччи многообразия, недавно был предложен в работе *Joulin A., Ollivier Y.* Curvature, concentration and error estimates for Markov chain Monte Carlo // Ann. Prob. 2010. V. 38. № 6. P. 2418–2442.

**Замечание.** Как уже неявно отмечалось в указании, полезно посмотреть на описанную модель с точки зрения моделей стохастической химической кинетики (см. замечание к задаче 19) с реакциями $1 \overset{\lambda}{\longrightarrow} 2$ и $2 \overset{\lambda}{\longrightarrow} 1$ и получить равновесную конфигурацию с помощью принципа максимума энтропии. Также на этом примере становится особенно прозрачной связь энтропии *Больцмана (функция Ляпунова прошкалированной кинетической динамики)* и энтропии *Санова (функционал действия в неравенствах больших уклонений для инвариантной меры)* – см. также замечание к задаче 19. На примере этой модели можно говорить о том, что в макросистемах возврат к неравновесным макросостояниям вполне допустим, но происходить это может только через очень большое время (*циклы Пуанкаре*), так что нам может не хватить отведенного времени, чтобы это заметить (*парадокс Цермело*). Напомним, что описанный выше случайный процесс обратим во времени. Однако наблюдается необратимая динамика относительной разности числа блох на собаках (*парадокс Лошмидта*). Но в таком случае можно удивляться также и тому, что газ, собранный в начальный момент в одной половине сосуда, с течением времени равномерно распределится по сосуду. Траектории будут обратимыми, но, поскольку, изначально система выведена из равновесия, то с большой вероятностью (тем большей, чем больше $N$) система будет стремиться вернуться в равновесие, где она стремится пребывать. Отсюда и возникает "стрела времени", см., например, концовку монографии *Опойцев В.И.* Нелинейная системостатика. М.: Наука, 1986.

**Задача 2 (ветвящийся процесс).** В колонию зайцев внесли зайца с необычным геном. Обозначим через $p_k$ – вероятность того, что в потомстве этого зайца ровно $k$ зайчат унаследуют этот ген ($k = 0, 1, 2, ...$). Это же распределение вероятностей характеризует всех последующих потомков, унаследовавших необычный ген. Будем считать, что каждый заяц дает потомство один раз в жизни в возрасте одного года (как раз в этом возрасте находился самый первый заяц с необычным геном в момент попадания в колонию).

Обозначим через $G(z)$ – производящую функцию распределения $p_k$, $k = 0, 1, 2, ...$, т.е. $G(z) = \sum_{k=0}^{\infty} p_k z^k$. Пусть $X_n$ – количество зайцев в возрасте одного года с необычным геном спустя $n$ лет после попадания в колонию первого такого зайца. Производящую функцию с.в. $X_n$ обозначим $\Pi_n(z) = E(z^{X_n})$.

**а)** Получите уравнение, связывающее $\Pi_{n+1}(z)$ с $\Pi_n(z)$ посредством $G(z)$.

**Указание.** Покажите, что $E(z^{X_{n+1}} | X_n) = [G(z)]^{X_n}$. Затем возьмите математическое ожидание от обеих частей равенства.

**б)** Покажите, что вероятность вырождения гена $q_n = P(X_k = 0; k \geq n) = \Pi_n(0)$. Существует ли предел $q = \lim_{n \to \infty} q_n$? Если существует, то найдите его.

**Указание.** Легко видеть, что функция $G(z)$ – выпуклая. Уравнение $z = G(z)$ имеет два корня: один в любом случае равен 1, другой $q \leq 1$. Если $\nu = G'(1) > 1$, то $q < 1$. Если $\nu \leq 1$, то $q = 1$.

**Замечание.** См. *Севастьянов Б.А.* Ветвящиеся процессы (серия "Теория вероятностей и математическая статистика"). М.: Наука, 1971; *Калинкин А.В.* Марковские ветвящиеся процессы с взаимодействием // УМН. 2002. Т. 57:2(344). С. 23–84. Частным, но важным случаем ветвящихся процессов, изученных в этой статье, являются модели стохастической химической кинетики, которые нам далее часто будут встречаться на протяжении всего раздела, в особенности, в замечании к задаче 19. В частности, представляется интересным (и не очень сложным) обобщить утверждение п. г) этого замечания на не много более общие модели из статьи А.В. Калинкина.

**Задача 3 (модель Шлёгля и теорема Т. Куртца).** Пусть $Y_n(t) \in \mathbb{Z}_+$, $t \geq 0$ – случайный процесс рождения и гибели с интенсивностями рождения и гибели ($h \to 0+$):[1]

$$P(Y_n(t+h) = j+1 | Y_n(t) = j) = n\lambda_n\left(\frac{j}{n}\right)h + o(h),$$

$$P(Y_n(t+h) = j-1 | Y_n(t) = j) = n\mu_n\left(\frac{j}{n}\right)h + o(h),$$

где

$$\lambda_n(x) = 1 + 3x\left(x - \frac{1}{n}\right), \quad \mu_n(x) = 3x + x\left(x - \frac{1}{n}\right)\left(x - \frac{2}{n}\right).$$

Введем прошкалированный процесс

$$X_n(t) = n^{1/4}\left(n^{-1}Y_n(n^{1/2}t) - 1\right), \, t \geq 0.$$

**а)** Положим $E_n = \left\{n^{1/4}\left(n^{-1}y - 1\right) : y \in \mathbb{Z}_+\right\}$, и определим (марковскую) полугруппу на банаховом пространстве (с нормой $\|\cdot\|$) достаточно хороших функций $f(x)$ над $E_n$ ($x \in E_n$):

$$T_n(t)f(x) \equiv E\left[f(X_n(t)) | X_n(0) = x\right].$$

Определим генератор полугруппы (линейный оператор), как (сходимость по норме $\|\cdot\|$)

---

[1] Этот процесс представляет собой модель (Шлёгля) стохастической химической кинетики:

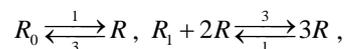

где $Y_n(t)$ – количество молекул типа $R$ в момент времени $t$. Над (под) стрелками написаны константы соответствующих реакций (см. замечание к задаче 19).

$$G_n f = \lim_{t \to 0+} \frac{1}{t}\{T_n(t)f - f\}.$$

Покажите, что

$$G_n f(x) = n^{3/2} \lambda_n (1 + n^{-1/4} x)\{f(x + n^{-3/4}) - f(x)\} + n^{3/2} \mu_n (1 + n^{-1/4} x)\{f(x - n^{-3/4}) - f(x)\}.$$

Положим $Gf(x) = 4f''(x) - x^3 f(x)$. Покажите, что этому генератору соответствует диффузионный процесс (отметим, что имеет место и единственность)

$$dX(t) = -X^3(t)dt + 2\sqrt{2}dW(t)$$

или (если известно $X(0)$)

$$X(t) = X(0) + 2\sqrt{2}W(t) - \int_0^t X^3(s)ds,$$

где $W(t)$ – винеровский процесс (напомним, в частности, что $W(t) \in N(0,t)$). Покажите, что для[2] $f \in C^2(\mathbb{R}) \cap C_c^1(\mathbb{R})$

$$\lim_{n \to \infty} \sup_{x \in E_n} |G_n f(x) - Gf(x)| = 0.$$

Теорема Т. Куртца утверждает, что если множество функций $f$, для которых это соотношение выполняется, достаточно богато (в нашем случае это так), и $X_n(0)$ слабо сходится к $X(0)$, то имеет место слабая сходимость соответствующих процессов:[3] $X_n \Rightarrow X$.

**б)** Покажите, что имеет место следующее представление в виде пуассоновских процессов со случайной заменой времени:

---

[2] Пространство дважды гладких функций с финитным носителем производной (производная отлична от нуля на ограниченном множестве).

[3] Введем $E$ полное сепарабельное метрическое пространство. Введем $D_E(\mathbb{R}_+)$ – пространство функций на $\mathbb{R}_+$ со значениями в $E$, непрерывных справа и имеющих пределы слева. Наделим пространство $D_E(\mathbb{R}_+)$ топологией Скорохода, при которой это пространство также будет полным сепарабельным метрическим пространством. Определим $\bar{C}(D_E(\mathbb{R}_+))$ – пространство ограниченных непрерывных функционалов над $D_E(\mathbb{R}_+)$. Будем говорить, что последовательность случайных процессов (со значениями в $E$, заданных на $\mathbb{R}_+$ и имеющих почти все траектории из $D_E(\mathbb{R}_+)$) $X_n$ слабо сходится к случайному процессу $X$ (краткая запись $X_n \Rightarrow X$), если $\lim_{n \to \infty} E[f(X_n)] = E[f(X)]$ для любого функционала $f \in \bar{C}(D_E(\mathbb{R}_+))$.

$$Y_n(t) = Y_n(0) + K_+\left(n\int_0^t \lambda_n\left(n^{-1}Y_n(s)\right)ds\right) - K_-\left(n\int_0^t \mu_n\left(n^{-1}Y_n(s)\right)ds\right),$$

$K_+$, $K_+$ – независимые пуассоновские процессы с одинаковым значением параметра равным 1 (напомним, в частности, что $K_+(t) \in Po(t)$). Покажите, что

$$X_n(t) = X_n(0) + n^{-3/4}\tilde{K}_+\left(n^{3/2}\int_0^t \lambda_n\left(1 + n^{-1}X_n(s)\right)ds\right) -$$

$$-n^{-3/4}\tilde{K}_-\left(n^{3/2}\int_0^t \mu_n\left(1 + n^{-1}X_n(s)\right)ds\right) + n^{3/4}\int_0^t \left(\lambda_n\left(1 + n^{-1}X_n(s)\right) - \mu_n\left(1 + n^{-1}X_n(s)\right)\right)ds,$$

где $\tilde{K}_+(u) = K_+ - u$, $\tilde{K}_-(u) = K_- - u$. Используя то, что

$$\left(n^{-3/4}\tilde{K}_+\left(n^{3/2}\cdot\right), n^{-3/4}\tilde{K}_-\left(n^{3/2}\cdot\right)\right) \Rightarrow (W_+, W_-),$$

где $(W_+, W_-)$ – независимые винеровские процессы, покажите, что $X_n \Rightarrow X$, где

$$X(t) = X(0) + W_+(4t) + W_-(4t) - \int_0^t X^3(s)ds.$$

Покажите, что такое представление предельного процесса совпадает (по распределению) с представление из предыдущего пункта (в общем случае это не всегда бывает так).

**Замечание.** Детали имеются в книге *Ethier N.S., Kurtz T.G.* Markov processes. Wiley Series in Probability and Mathematical Statistics: Probability and Mathematical Statistics. John Wiley & Sons Inc., New York, 2005.[4] Хочется отметить, что результаты, описанные в этой книге, будут полезны для решения многих задач этого раздела. В рассмотренной задаче мы разобрали простейший пример, который демонстрирует основные конструкции. Все это можно существенно обобщить. В частности, хотелось бы отметить приложения описанной техники к изучению ветвящихся процессов (глава 9) и строгому выводу модели популяционной генетики Райта–Фишера (глава 10). Популяционная генетика и ее окрестности хорошо ложатся в тематику данного раздела, однако в виду достаточно большого количества материала мы лишь ограничимся ссылкой на монографию *Свирежев Ю.М., Пасеков В.П.* Основы математической генетики. М.: Наука, 1982 и популярную книгу *Горбань А.Н., Хлебопрос Р.Г.* Демон Дарвина. Идея оптимальности и естественный отбор. М.: Наука, 1988 (в электроном виде имеется второе издание: Красноярск, 1998). Также можно отметить недавно переведенную довольно оригинальную книгу *Хойл Ф.* Математика эволюции. М.–Ижевск: НИЦ "РХД", 2012.

---

[4] Помимо этой книги, также хотелось бы обратить внимание на книгу *Гардинер К.В.* Стохастические модели в естественных науках. М.: Мир, 1986. В этой книге изложены менее строгие, но более простые способы получения предельных уравнений для различных марковских динамик. В частности, описаны разложение Крамерса–Мойала и разложение ВанКампена (по обратному размеру системы).

**Задача 4 (Кинетика социального неравенства и предельные формы). а)**\*\* В некотором городе живет $N \gg 1$ жителей (четное число). В начальный момент у каждого жителя имеется по $\bar{s}$ монеток. Каждый день жители случайно разбиваются на пары. В каждой паре жители скидываются по монетке (если один или оба участника банкроты, то банкрот не скидывается, в то время, как не банкрот, в любом случае, обязан скинуть монетку). Далее в каждой паре случайно разыгрывается победитель, который и забирает "призовой фонд". Обозначим через $c_s(t)$ – долю жителей города, у которых ровно $s$ ($s = 0,...,\bar{s}N$) монеток на $t$-й день. Покажите, что

$$\exists\ a > 0:\ \forall\ \sigma > 0, t \geq a(\bar{s})\ln N \rightarrow P\left(\|c(t) - c^*\|_2 \geq \frac{2\sqrt{2} + 4\sqrt{\ln(\sigma^{-1})}}{\sqrt{N}}\right) \leq \sigma,$$

$$\exists\ b, D > 0:\ \forall\ \sigma > 0, t \geq b(\bar{s})\ln N \rightarrow P\left(\|c(t) - c^*\|_1 \geq D\sqrt{\frac{\ln^2 N + \ln(\sigma^{-1})}{N}}\right) \leq \sigma,$$

где $c_s^* \simeq C\exp(-s/\bar{s})$, а $C \simeq 1/\bar{s}$ находится из условия нормировки $\sum_{s=0}^{\bar{s}N} C\exp(-s/\bar{s}) = 1$. Таким образом, кривая (предельная форма) $C\exp(-s/\bar{s})$ характеризует распределение населения по богатству на больших временах.

**Указание.** Для решения этой задачи полезно рассмотреть схожий процесс, в котором каждой паре жителей приписан свой (независимый) пуассоновский будильник (звонки происходят в случайные моменты времени, соответствующие скачкам пуассоновского процесса; параметр интенсивность этого пуассоновского процесса называют интенсивностью/параметром будильника). Все будильники "приготовлены" одинаково: у всех у них одна и та же интенсивность $\lambda N^{-1}$. Далее следует погрузить задачу в модель стохастической химической кинетики с бинарными реакциями и воспользоваться результатом из замечания к задаче 19. Наиболее технически сложными моментами в получении указанного в условии задачи результата является оценка mixing time $\sim \ln N$ и получение поправки под корнем $\ln^2 N$. Насколько нам известно, пока это еще нигде аккуратно не обосновано.

**Замечание.** Название модели мы взяли из одноименной статьи К.Ю. Богданова[5] в журнале "Квант". В этой статье предлагаются и другие правила обмена. О возможных обобщениях этой модели также можно посмотреть в работах *Dragulescu A., Yakovenko*

---
[5] К.Ю. Богданов выпустил очень познавательную брошюрку "Прогулки с физикой" в серии "Библиотечка Квант", в которую вошел этот сюжет и многие другие. Также у К.Ю. Богданова мы позаимствовали идею одной из стохастических динамик, приводящих к модели хищник–жертва и системе уравнений Лотки–Вольтера, см. задачу 10.

*V.M.* Statistical mechanics of money // The European Physical Journal B. 2000. V. 17. P. 723–729 и в других работах этих авторов, а также в работах A.M. Chebotarev'a.

**б)\* (задача Булгакова–Маслова о разбрасывании червонцев в варьете)** В некотором городе живет $N \gg 1$ жителей (изначально банкротов). Каждый день одному из жителей (случайно выбранному в этот день) дается одна монетка. В обозначениях п. а), определите предельную форму для $\{c_s(\overline{s}N)\}_{s=0}^{\overline{s}N}$. Сколько надо случайно раздать монеток, чтобы с вероятностью $\geq 1-\sigma$ в городе бы не было банкротов?[6]

**Указание.** Покажите, что распределение $\{c_s(\overline{s}N)\}_{s=0}^{\overline{s}N}$ – мультиномиальное. Более того, оно в точности совпадает со стационарным распределением марковского процесса из предыдущего пункта. В работах В.П. Маслова (в основном в журнале Мат. заметки) собрано большое число тонких результатов о концентрации различных (урновых) мер, в том числе возникающих в приложениях. Такие меры, как правило, имеют комбинаторную природу, и описывают способы распределения шаров по урнам при различного рода линейных ограничениях. В нашем случае урны отвечают $s = 0,...,\overline{s}N$. А ограничения – закон сохранения числа жителей и числа монеток (задают множество (inv) в замечании к задаче 19)

$$\sum_{s=0}^{\overline{s}N} n_s(t) = N, \ \sum_{s=0}^{\overline{s}N} s n_s(t) = \overline{s}N,$$

где $n_s(t)$ – количество жителей, у которых в ровно $s$ монеток на $t$-й день.

**в)\*\* (закон Ципфа–Парето и процесс Юла)** В некотором городе живет неограниченно много жителей (изначально банкротов), которые могут участвовать в "освоении" монеток. База индукции: сначала выбирается один житель, он получает одну монетку. Шаг индукции: в $(k+1)$-й день выбирается очередной новый житель (отличный от $k$ уже выбранных), он получает одну монетку, вторая монетка с вероятностью $\alpha < 1$ равновероятно отдается одному из $k$ старых жителей, а с вероятностью $1-\alpha$ эта монетка отдается одному из $k$ старых жителей с вероятностью, пропорциональной тому, сколько у него уже есть монеток, т.е. по принципу "деньги к деньгам" (в моделях роста Интернета этот принцип называют "preferential attachment"). Поясним примером. Пусть таких старых жителя три ($k = 3$). У первого 3 монетки у второго 1 и у третьего 1. Тогда с вероятностью 3/5 монетка попадет к первому, с вероятностью 1/5 ко второму и 1/5 к третьему. Аналогично и в общем случае. В обозначениях п. а), определите предельную форму для $\{c_s(t)\}_{s=1}^{t}$, $t \gg 1$.

**Указание.** Изучать стохастическую динамику в этой задаче заметно сложнее, чем в двух предыдущих. Поэтому обычно работают в приближении среднего поля. В данном контексте это означает, что $n_s(t) \simeq E[n_s(t)]$. Далее выписывают на $n_s(t)$ систему зацепляющихся обыкновенных дифференциальных уравнений. Автомодельное притягивающее решение этой системы ищут в виде $n_s(t) \sim c_s^* t$. После разрешения соответствующих уравнений получают степенной закон для зависимости $c_s^* \sim s^{-(3+\alpha/(1-\alpha))}$. Такой подход применительно к моделям роста интернета (и изучения степенного закона для распределения степеней вершин) довольно часто сейчас встречается (в том числе и в учебной литературе). В частности, этот подход описан в обзоре *Mitzenmacher M.* A Brief History of Generative Models for Power Law and Lognormal Distributions // Internet Mathematics. 2003. V. 1. no. 2. P. 226–251, и рассматривается в книге Hopkroft–Kannan'a (Computer Science Theory for the Information Age – не вышедшая пока книга, свободно доступная в интернете). По-сути, в этом пункте нами был описан процесс, возникающий в работах 20-х годов XX века по популяционной генетике, получивший названия процесса Юла (см., например, обзор *Newman M.E.O.* Power laws, Pareto distributions and Zipf's law // Contemporary physics. 2005. V. 46. no. 5. P. 323–351).

---
[6] Оказывается, ответ на этот вопрос помогает организовывать (определять сколько нужно выбрать точек старта) метод случайного мультистарта в глобальной оптимизации, см., например, *Жиглявский А.А., Жилинскас А.Г.* Методы поиска глобального экстремума. М.: Наука, 1991.

**Замечание.** Описанные модели восходят к работе конца XIX века В. Парето, в которой была предпринята попытка объяснить социальное неравенство и к работе Г. Ципфа конца 40-х годов XX века, в которой была отмечена важность степенных законов в "Природе". Эти законы для большей популярности иногда преподносят, как принцип Парето или принцип 80/20 (80 % результатов проистекают всего лишь из двадцати процентов причин) – такие пропорции отвечают $c_s^* \sim s^{-2.1}$ (см. Newman M.E.O.). Приведем примеры (не совсем, правда, точные): 80 % научных результатов получили 20 % ученых, 80 % пива выпило 20 % людей и т.п. Сейчас много исследований во всем мире посвящено изучению возникновения в самых разных приложениях степенных законов (распределение городов по населению, коммерческих компаний по капитализации, автомобильных пробок по длине). Особенно бурный рост возник в связи с изучением роста больших сетей (экономических, социальных, интернета), см., например, книги и работы M.O. Jackson'a. В России это направление также представлено, например, в Яндексе в отделе А.М. Райгородского (Гречников–Остроумова–Рябченко–Самосват, Леонидов–Мусатов–Савватеев), в ИПМ РАН в группе А.В. Подлазова (собственно, модель п. в) была взята из статьи Подлазова А.В. в сборнике Новое в синергетике. Нелинейность в современном естествознании. М.: ЛИБРОКОМ, 2009. С. 229–256).

**Задача 5 (обезьянка и печатная машинка; закон Ципфа–Мандельброта).** На печатной машинке $n+1$ символ, один из символов пробел. Обезьянка на каждом шаге случайно (независимо и равновероятно) нажимает один из символов. Прожив долгую жизнь, обезьянка сгенерировала текст огромной длинны. По этому тексту составили словарь. Этот словарь упорядочили по частоте встречаемости слова (слова – это всевозможные набор букв без пробелов, которые хотя бы раз встречались в тексте обезьянки между какими-то двумя пробелами). Так на первом месте в словаре поставили самое часто встречаемое слово, на второе поставили второе по частоте встречаемости и т.д. Номер слова в таком словаре называется рангом и обозначается буквой $r$. Покажите, что предельная форма кривой, описывающей распределение частот встречаемости слов от рангов, имеет вид

$$\text{Частота}(r) \simeq \frac{C}{(r+B)^{\alpha}},$$

где

$$\alpha = \frac{\log(n+1)}{\log(n)},\ B = \frac{n}{n-1},\ C = \frac{n^{\alpha-1}}{(n-1)^{\alpha}}.$$

**Замечание.** Такой вывод закона Ципфа (с поправкой) был одним из двух, предложенных Б. Мандельбротом (создателем науки о фракталах[7]). Этот вывод можно найти во многих источниках. Благодаря Википедии, наиболее популярной на эту тему оказалась статья *Wentian Li* Random Text Exibit Zipf's–Law–Like Word Frequency Distribution // IEEE Transactions of Information Theory. 1992. V. 38(6). P. 1842–1845. В работе (см. также Conrad–Mitzenmacher) *Бочкарева В.В., Лернера Э.Ю.* Закон Ципфа для случайных текстов с неравными вероятностями букв и пирамида Паскаля // Известия вузов. Математика. 2012. № 12. С. 30–33 приводится обобщение этой модели на случай неравных вероятностей букв. В другой работе этих авторов рассматривается обобщение, связанное с рассмотрением ситуации, когда последовательность нажимаемых символов не i.i.d., а образует цепь Маркова (впрочем, такого рода обобщения довольно популярны, и рассматривались рядом других авторов). В работе *Шрейдер Ю.А.* О возможности теоретического вывода статистических закономерностей текста (к обоснованию закона Ципфа) // Проблемы передачи информации. 1967. Т. 3. № 1. С. 57–63 была предложена довольно общая схема, приводящая к закону Ципфа, в которую можно погрузить[8] и обезьянку с печатной машинкой. А именно, предположим, что динамика (порождения слов в большом тексте) такова, что вероятность того, что в тексте из $N \gg 1$ слов $x_1$ (первое по порядку слово в ранговом словаре) встречалось $N_1$ раз, $x_2$ (второе по порядку слово в ранговом словаре) встречалось $N_2$ раза и т.д. есть

$$\sim \frac{N!}{N_1! N_2! \ldots} \exp\left(-\eta \sum_{k \in \mathbb{N}} N_k E_k\right),$$

где $E_k$ – число букв в слове с рангом $k$. Часто считают, что $\eta = 0$, но зато динамика такова, что число слов и число букв становятся асимптотически (по размеру текста) связанными (закон больших чисел). Таким образом, к закону сохранения $\sum_{k \in \mathbb{N}} N_k = N$ добавляется приближенный закон сохранения $\sum_{k \in \mathbb{N}} N_k E_k \simeq \bar{E} N$ ($\bar{E}$ – среднее число букв в слове). Поиск предельной формы приводит к задаче (воспользовались формулой Стирлинга и методом множителей Лагранжа)[9]

$$\sum_{k \in \mathbb{N}} \{N_k \ln N_k + \lambda E_k N_k\} \to \min_{\substack{N_k \geq 0 \\ \sum_{k \in \mathbb{N}} N_k = N}},$$

где $\lambda$ – либо равняется $\eta$, либо является множителем Лагранжа к ограничению $\sum_{k \in \mathbb{N}} N_k E_k \simeq \bar{E} N$. Стоит отметить, что к аналогичной задаче (с $E_k = k$) приводит поиск предельной формы в модели "Кинетика социального неравенства" из п. а) задачи 4. Решение нашей задачи дает $N_k = \exp(-\mu - \lambda E_k)$, где $\mu$ – множитель Лагранжа к ограничению $\sum_{k \in \mathbb{N}} N_k = N$. Далее, считают, что $r(E)$ число различных используемых слов с числом букв не большим $E$ приближенно представимо в виде $r(E) \simeq a^E$. Тогда $N_k \sim k^{-\gamma}$, где $\gamma = \lambda / \ln a$.

---

[7] Интерес Б. Мандельброта к степенным законам (законам типа Ципфа) был закономерен, поскольку степенные законы как раз проявляют свойства масштабной инвариантности (самоподобности), свойственной фракталам, и часто возникают при описании различных критических явлений (см., например, обзор M.E.O. Newman'a, 2005). В популярных книгах Б. Мандельброта, большинство из которых переведены на русский язык, также периодически встречаются сюжеты на эту тему. Вообще, здесь бы хотелось отметить, что очень многие сложные сети, возникающие в природе (капиллярная система, ветви деревьев и т.п.) имеют в своей основе определенные закономерности, формирующиеся при росте (напомним, что социальные сети и интернет моделируются с помощью модели случайного роста "предпочтительного присоединения"). В результате возникающая конфигурация находится из некоторого вариационного принципа, приводящего, также как и рассмотренной задаче к степенным законам. Яркий пример имеется, например, в статье *West G.B., Brown J.H. Enquist B.J.* A General Model for the Origin of Allometric Scaling Laws in Biology // SCIENCE. 1997. V. 276. P. 122–126. Представляется интересным искать аналоги "законов физики" при изучения роста не биологических сетей. То есть небольшой набор универсальных правил, которые в различных сочетаниях и в различных контекстах дают многообразие возникающих в приложениях сетей больших размеров, и отражают их основные (статистические) свойства. Один такой закон (предпочтительного присоединения) уже открыт (у этого закона много вариаций, за счет которых и удается подогнать ту или иную модель под экспериментальные результаты). Работа в этом направлении начата относительно недавно (~15 лет назад), и, кажется, что здесь все самое интересное впереди.

[8] Это полезно показать. Также представляется полезным изучение других примеров (в том числе упомянутых) на предмет возможности их погружения в эту схему.

Некоторые обобщения закона Ципфа–Мандельброта с обсуждением также и лингвистической стороны дела имеются в работе В.П. Маслова с дочерью. Интересным представляется связь, недавно обнаруженная Ю.И. Маниным, между законом Ципфа и колмогоровской сложностью. Отметим, что ранее на такую связь также было указано В.В. Вьюгиным.

**Задача 6 (Heaps' law)\*.** Предположим, что сгенерирован достаточно большой текст, в котором слова случайно (независимо) выбирались из словаря с вероятностями, зависящими от ранга слова по закону Ципфа–Мандельброта (см. предыдущую задачу). Покажите, что если сгенерированный текст имеет всего $n$ слов, то число различных слов $m$ в нем будет $\sim n^{1/\alpha}$. Постарайтесь получить более точное описание с.в. $m$.

---

[9] Интересно было бы также исследовать и концентрацию меры вокруг этой предельной формы.

**Замечание.** См., например, *Grootjen F.A., Leijenhorst D.C., van der Weide Th.P.* A formal deviation of Heaps' law // Inform. Sciences. 2005. V. 170(24). P. 263–272.

**Задача 7 (распределение букв по частоте встречаемости; С.М. Гусейн-Заде).** Предположим, что есть некоторая эргодическая динамика (детерминированная) на единичном симплексе с равномерной инвариантной мерой, описывающая эволюцию частот встречаемости букв в текстах. Эту динамику можно понимать, как эволюцию состояния некоторого лингвистического мира. Мы считаем, что это динамка развивается в медленном времени (годы). В то же время, в быстром времени (дни), при заданном "состоянии мира" (частотах встречаемости букв) случайно (согласно этим частотам) генерируется (с постоянной интенсивностью) большое количество текстов. За много лет накопилось огромное количество таких текстов. Их объединили все вместе в один огромный текст и посчитали частоты встречаемости различных букв в этом объединенном тексте. Считая, что всего имеется $n$ различных букв, покажите, что $r$-я по величине частота встречаемости приблизительно равна (ранговое распределение частот букв)

$$\frac{1}{n}\sum_{k=r}^{n}\frac{1}{k} \approx \frac{1}{n}\left(\ln(n+1) - \ln r\right).$$

**Замечание.** Впрочем, описанный закон в определенном смысле ничем не примечателен по сравнению с рядом известных конкурирующих законов, см. *Гельфанд М.С., Минь Чжао* О ранговых распределениях частот букв в естественных языках // Пробл. передачи информ. 1996. Т. 32. № 2. С. 89–95.

**Задача 8 (случайный рост диаграмм, предельные формы).** Под диаграммой будем понимать горку (на плоскости) составленную из одинаковых квадратных кирпичей. Слово "горка" означает, что ряд из кирпичей, поставленных друг на друга, не ниже ряда из кирпичей, примыкающего справа.

**а)\*\* (диаграммы Ричардсона)** Пусть горка растет по принципу: новый кирпич кладется с равной вероятностью на одну из допустимых позиций, где под допустимой позицией понимается "уголок" либо две "крайние" позиции (над самым левым верхнем кирпичом, справа от самого правого нижнего кирпича – далее мы будем называть эти позиции, соответственно, верхний и нижний уголок). Если положить новый кирпич в такой уголок на горке, то он будет касаться горки двумя своими гранями (левой и нижней). Пусть случайно, как описано выше, положили $n \gg 1$ кирпичей. Перейдя в систему координат с прошкалированными фактором $\sqrt{n}$ осями и разместив левый нижний кирпич горки в начале координат, покажите, что в пределе при $n \to \infty$ форма горки (предельная форма) будет неслучайной кривой, описываемой уравнением $\sqrt{x} + \sqrt{y} = \sqrt[4]{6}$, где константа $\sqrt[4]{6}$ была выбрана из условий нормировки площади под кривой на 1.

**Замечание.** См. статью *Rost H.* Nonequilibrium behavior of a many particle process: Density profile and local equilibria // Probability Theory and Related Fields. 1981. V. 58. no. 1. P. 41–53.

**б)*** Как надо правильно прошкалировать оси и какая при этом получится предельная форма, если новый кирпич кладется с вероятностью $\alpha$ в нижний уголок и с вероятностью $1-\alpha$ с равновероятно опускается на один из рядов кирпичей, поставленных друг на друга, и далее скатывается в уголок, находящийся на этом уровне? Как изменится ответ, если новый кирпич равновероятно опускается на один из рядов кирпичей, поставленных друг на друга, или нижний уголок, и далее скатывается в уголок, находящийся на этом уровне?

**Замечание.** Следует сравнить с п. в) задачи 4. Предельные формы различных диаграмм активно изучались и изучаются в работах А.М. Вершика и его учеников. Часто эти диаграммы имеют вполне содержательную интерпретацию (например, диаграммы Юнга). Мера на диаграммах часто задается без описания динамики, т.е. не описывается рост диаграмм. Таким примером является изучение предельной формы диаграмм Юнга с равномерной мерой (рассматривают и другие меры, например, Планшереля). Сюда же можно отнести и статистику выпуклых ломанных (Арнольд–Вершик–Синай). Примечательно, что в обоих случаях предельная форма определяется из решения вариационной задачи, которая получается исходя из максимизации (в пространстве кривых) должным образом прошкалированного логарифма числа диаграмм (выпуклых ломанных), находящихся в малой окрестности заданной кривой (аналог действия/энтропии в теоремах о больших уклонениях типа Санова[10]). Причем ответ не зависит от того, как выбирается это окрестность (в предположении ее малости), важно лишь, чтобы окрестности строились по одним и тем же правилам для разных кривых.

**Задача 9 (распределенные модели стохастической химической кинетики)*.** На окружности единичного периметра на расстоянии $\varepsilon$ друг от друга расположены дома. В доме каждом доме живет по одной собаке и одной кошке. Для определенности, в начальный момент у каждой собаки $N$ блох, а на кошках блох в начальный момент нет. С интенсивностью $\lambda$ (см. задачу 1 "парадокс Эренфестов") каждая блоха независимо от остальных перескакивает с собаки (кошки) на соседа по дому, т.е. кошку (собаку), а с вероятностью $\lambda_{-1}/\varepsilon + \mu/\varepsilon^2$ ($\lambda_{+1}/\varepsilon + \mu/\varepsilon^2$) на собаку (кошку) из соседнего дома с на единицу меньшим (большим) номером. Покажите (используя теорему Куртца), что при $N \to \infty$, $\varepsilon \to 0+$ п.н. существует (детерминированный) предел (причем, на каждом конечном отрезке времени можно говорить о слабой сходимости случайных процессов)

$$c^{dog}(t,x) = \lim_{\substack{N\to\infty,\,\varepsilon\to 0+ \\ x/\varepsilon \le k < x/\varepsilon+1}} \frac{n_k^{dog}(t)}{N},$$

---

[10] Для определения предельной формы диаграмм Юнга мы получаем в точности задачу энтропийно-линейного программирования в функциональном пространстве, следует сравнить с замечанием к задаче 19.

где $n_k^{dog}(t)$ – количество блох у собаки в доме с номером $k$ в момент времени $t$, аналогично, через $n_k^{cat}(t)$, определяется $c^{cat}(t,x)$. Покажите, что введенные плотности удовлетворяют следующей системе УЧП

$$\frac{\partial c^{dog}}{\partial t} = \mu \frac{\partial^2 c^{dog}}{\partial x^2} + (\lambda_{-1} - \lambda_{+1})\frac{\partial c^{dog}}{\partial x} + \lambda \cdot (c^{cat} - c^{dog}),$$

$$\frac{\partial c^{cat}}{\partial t} = \mu \frac{\partial^2 c^{cat}}{\partial x^2} + (\lambda_{-1} - \lambda_{+1})\frac{\partial c^{cat}}{\partial x} + \lambda \cdot (c^{dog} - c^{cat}),$$

унаследовавшей линейный закон сохранения от стохастической динамики

$$\int_0^1 \{c^{dog}(t,x) + c^{cat}(t,x)\} dx \equiv 1. \qquad \text{(inv)}$$

Покажите, что стационарное распределение описанного марковского процесса с носителем (inv) имеет асимптотическое представление

$$\sim \exp\left(-\frac{N}{\varepsilon} \int_0^1 \{c^{dog}(t,x)\ln(c^{dog}(t,x)) + c^{cat}(t,x)\ln(c^{cat}(t,x))\} dx\right),$$

Причем функционал (пространственная энтропия), стоящий при $N/\varepsilon$, является функционалом Ляпунова выписанной ранее системы УЧП. Исследуйте поведение системы на больших временах.

**Замечание.** Для решения этой задачи (и ряда других задач этого раздела) можно рекомендовать монографии

*Ethier N.S., Kurtz T.G.* Markov processes. Wiley Series in Probability and Mathematical Statistics: Probability and Mathematical Statistics. John Wiley & Sons Inc., New York, 2005.

*Гардинер К.В.* Стохастические модели в естественных науках. М.: Мир, 1986 (главы 7, 8).

В данной задаче рассмотрен, пожалуй, простейший пример пространственно распределенной системы. Такие системы часто возникают, например, в математической биологии при описании эволюции пространственно распределенных взаимодействующих видов. Причем, как правило, за счет более сложного взаимодействия и(или) неоднородной (анизотропной) диффузии и(или) сноса (в том числе, со скоростью/интенсивностью, зависящей от концентраций видов), при скейлинге получаются уже нелинейные эволюционные уравнения (параболического типа, если есть диффузия). При этом равновесная конфигурация, возникающая на больших временах часто может быть также проинтерпретирована в терминах равновесия по Нэшу и Дарвиновского отбора (не распределенный вариант, поясняющий отмеченные связи, имеется в задачах 10, 15, 16). К этому направлению можно отнести и "принцип эволюционной оптимальности", активно развиваемый в работах В.Н. Разжевайкина.

**Задача 10 (модель хищник–жертва)**. В некотором сказочном лесу, имеющем форму единичного тора (прямое произведение двух окружностей единичного периметра), который можно понимать как единичный квадрат с отождествленными противоположными сторонами, проведена решетка с квадратными ячейками размера $\varepsilon \times \varepsilon$. В узлах этой решетки в начальный момент времени случайно (независимо и равновероятно) распределены $n_з = \kappa_з \varepsilon^{-1}$ зайцев и $n_в = \kappa_в \varepsilon^{-1}$ волков. Волки и зайцы начинают независимо случайно блуждать по решетке: с интенсивностью $\lambda_з \varepsilon^{-1}$ заяц покидает текущий узел решетки, выбирая для перехода равновероятно один из четырех соседних узлов, аналогично (с $\lambda_в \varepsilon^{-1}$) ведут себя волки. Кроме того, с интенсивностью $\mu_з^+$ заяц создает в текущем узле потомка, который сразу же начинает жить независимой жизнью по общим правилам. А волки с интенсивностью $\mu_в^-$ выбывают из системы. Если в результате блуждания волк и заяц окажутся в одном узле, то мгновенно происходит реакция: волк съедает зайца и производит потомка, который сразу же начинает жить независимой жизнью по общим правилам. В каком смысле можно говорит о том, что в пределе при $\varepsilon \to 0+$ (на конечных отрезках времени) динамика изменения общих численностей волков и зайцев соответствует модели стохастической химической кинетики (см. замечание к задаче 19)

$$[З] \xrightarrow{\mu_з^+} 2[З],\ [B]+[З] \xrightarrow{K} 2[B],\ [B] \xrightarrow{\mu_в^-} 0,$$

где константа реакции $K$ определяется по введенным выше параметрам? Покажите, что закон действующих масс приводит к системе Гульдберга–Вааге (получающейся при каноническом скейлинге, см. замечание к задаче 19)

$$\frac{dc_з}{dt} = \mu_з^+ c_з - K c_в c_з,$$

$$\frac{dc_в}{dt} = K c_в c_з - \mu_в^- c_в.$$

Эту систему уравнений принято называть системой (Лотки–)Вольтера.

**Замечание.** Несложно показать, что в стохастической модели (большое, но конечное, число волков и зайцев) с вероятностью 1 на больших временах система "свалится" в поглощающее состояние (без волков). При этом число зайцев в этом состоянии либо будет равняться нулю, либо стремиться к бесконечности со временем. При этом система Лотки–Вольтера имеет нетривиальное положение равновесия $\left(c_з^*, c_в^*\right) = \left(\mu_в^-/K, \mu_з^+/K\right)$ являющееся центром (в окрестности которого система колеблется с частотой 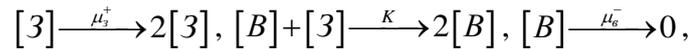 $\sqrt{\mu_з^+ \mu_в^-}$) причем не только в линейном приближении (локально), но и глобально. Это следует из консервативности системы. Первый интеграл имеет вид

$$\mu_в^- \ln c_з + \mu_з^+ \ln c_в - K \cdot (c_з + c_в) \equiv \mathrm{const}\left(c_з(0), c_в(0)\right).$$

Все это означает негрубость системы (чувствительность к возмущениям). Таким образом, для данной модели реакций стохастической химической кинетики мы имеем не перестановочность порядка взятия предельных переходов $\lim_{t\to\infty}\lim_{N\to\infty} \ne \lim_{N\to\infty}\lim_{t\to\infty}$, где число молекул (агентов) = волков + зайцев. Следует сравнить это замечание с замечанием к задаче 19, в котором приводится достаточное условие перестановочности. В данной модели можно лишь говорить о том, что система (Лотки–)Вольтера аппроксимирует соответствующую модель стохастической химической кинетики лишь на конечном промежутке времени. Этот промежуток тем больше, чем больше молекул (агентов). Более того, длина промежутка неограниченно возрастает с возрастанием числа агентов. Этот конкретный пример рассматривался в большом числе книг, переведенных на русский язык (Николис–Пригожин, К.В. Гардинер, В.Б. Занг).

Если рассматривать более общие биологические модели взаимодействия $n$ видов, то в пространственно однородном случае возникает аналогичная СОДУ

$$\frac{dc_i}{dt} = c_i f_i(c_1,...,c_n), \ i=1,...,n.$$

Функции $f_i(c_1,...,c_n)$ – называют мальтузианскими функциями. Если предположить, что имеется устойчивое положение равновесия $\bar{c} = \underbrace{(c_1,...,c_m,0,...,0)}_{n}$, то[11] необходимо $f_i(\bar{c}) = 0$ при $i=1,...,m$ и $f_i(\bar{c}) \le 0$ при $i=m+1,...,n$. Это означает (принцип эволюционный оптимальности), что

$$из \ \bar{c}_i > 0 \ следует \ f_i(\bar{c}) = \max_{j=1,...,n} f_j(\bar{c})$$

или, что тоже самое (задача дополнительности[12])

$$\bar{c}_i \cdot \left(\max_{j=1,...,n} f_j(\bar{c}) - f_i(\bar{c})\right) = 0.$$

Биологический смысл этого условия заключается в том, что выжившие в том или ином состоянии виды обязаны иметь максимальное значение мальтузианских коэффициентов

---

[11] Поскольку матрица Якоби имеет верхне-треугольный блочный вид

$$J = \left\|\frac{\partial(c_i f_i(c))}{\partial c_j}\right\|_{c=\bar{c}} = \begin{pmatrix} A & B \\ 0 & D \end{pmatrix},$$

спектр матрицы (по теореме Лапласа) задается спектром двух подматриц $A$ и $D = \text{diag}\{f_i(\bar{c})\}_{i=m+1}^{n}$. Для устойчивости необходима и достаточна внутренняя устойчивость (матрицы $A$), означающая устойчивость к изменению численностей присутствующих в равновесии видов, и внешняя (матрицы $D$), означающая устойчивость к привнесению в систему из вне видов, которые в равновесие не представлены. Последнее означает, что $f_i(\bar{c}) \le 0$ при $i=m+1,...,n$ (необходимое условие устойчивости матрицы $D$).

[12] Во многих важных приложениях это задача переписывается как (монотонное) вариационное неравенство, которое в случае потенциальности векторного поля $f$ переписывается как задача выпуклой оптимизации (см., например, задачу 15).

среди всех потенциально допустимых, вычисленных в этом состоянии. Именно эти коэффициенты характеризуют силу видов в ее дарвиновском понимании. Подробнее об этом и о возможных обобщениях можно посмотреть, например, в учебном пособии *Разжевайкин В.Н.* Анализ моделей динамики популяций. М.: МФТИ, 2010.

В популяционной теории игр рассмотренному выше придается немного другая интерпретация. А именно $c_i$ – концентрация (доля) игроков использующих стратегию $i$, $g_i(c_1,...,c_n)$ – выигрыш от использования стратегии $i$ при заданных пропорциях. Каждый игрок может выбирать любую из $n$ стратегий. Точнее, каждый игрок (агент), использующий стратегию $i$ ($i=1,...,n$), независимо от остальных с интенсивностью[13] $\lambda c_j \left[ g_j(c) - g_i(c) \right]_+$ меняет свою стратегию на $j \neq i$ (система представляет собой модель стохастической химической кинетики с унарными реакциями и с непостоянными коэффициентами реакций, см. задачу 15). Если число игроков стремится к бесконечности так, что корректно определены начальные концентрации $c(0)$, то концентрации корректно определены и при $t \geq 0$ (теорема Т. Куртца), причем вектор-функция $c(t)$ удовлетворяет выписанной ранее СОДУ с $f_i = g_i - \bar{g}$, где $\bar{g}(c) = \sum_{i=1}^{n} c_i g_i(c)$, которую (для заданного протокола) называют динамикой репликаторов. В такой интерпретации полученное равновесие приобретает смысл равновесия Нэша в популяционной игре. А принцип эволюционный оптимальности можно понимать как условие равновесия по Нэшу. Все необходимые детали имеются в книге *Sandholm W.* Population games and Evolutionary dynamics. Economic Learning and Social Evolution. MIT Press; Cambridge, 2010.

Предложенный в задаче скейлинг далеко не единственно возможный (см., например, задачу 9). Здесь мы приведем еще один вариант. Если считать, что $n_з = \kappa_з \varepsilon^{-2}$, $n_в = \kappa_в \varepsilon^{-2}$, и предположить, что находящиеся в одном узле волк и заяц прореагируют не мгновенно, а с интенсивностью $K$ (в остальном также как и раньше с поправкой, что $\lambda_з := \lambda_з \varepsilon^{-2}$, $\lambda_в := \lambda_в \varepsilon^{-2}$), то для любого $t \geq 0$ и $x \in [0,1)^2$ существует

$$c_з(t,x) = \lim_{\substack{\varepsilon \to 0+ \\ x_i/\varepsilon \leq k_i < x_i/\varepsilon+1, i=1,2}} E\left[ n_з(t;k_1,k_2) \right],$$

---

[13] Такой протокол (способ изменения своей стратегии) называют попарным сравнением. Но можно рассматривать и многие другие протоколы, см., например, книгу W. Sandholm'a и задачи 14, 15. Общей замечательной особенностью протоколов, отражающих рациональность игроков, является (при весьма общих условиях, но не всегда – см., например, задачу 14) одинаковость аттрактора, что можно грубо сформулировать так: "все дороги ведут в Рим". То есть при весьма общих условиях, если игроки стремятся менять, исходя из текущих наблюдаемых концентраций, свои стратегии в сторону увеличения собственного выигрыша, то соответствующая стохастическая динамика имеет инвариантную (стационарную) меру, которая, часто концентрируется, с ростом числа игроков, около (устойчивых) равновесий Нэша, множество которых не зависит от выбранного протокола. В свою очередь, множество, на котором концентрируется стационарная мера (в общем случае это могут быть и не равновесия Нэша, до которых динамика может "не дойти", скажем, "закружившись в цикл", которому соответствует "кратер вулкана" – см., например, задачу 14) можно рассматривать, как аттрактор соответствующей СОДУ.

где $n_з(t; k_1, k_2)$ – число зайцев в момент времени $t$ в узле решетки с номером $(k_1, k_2)$ (аналогично определяется $c_в(\tau, x)$). Причем эти функции удовлетворяют следующей системе (Лотки–Вольтера с диффузией) уравнений реакции-диффузии

$$\frac{\partial c_з}{\partial t} = \lambda_з \Delta c_з + \mu_з^+ c_з - K c_в c_з,$$

$$\frac{\partial c_в}{\partial t} = \lambda_в \Delta c_в + K c_в c_з - \mu_в^- c_в.$$

Детали и обобщения имеются в статье *Малышев В.А., Пирогов С.А.* Обратимость и необратимость в стохастической химической кинетике // Успехи мат. наук. 2008. Т. 63. вып. 1(379). С. 4–36. В частности, если жизнь "разворачивается" на $d$-мерном торе, то изменения будут только в определении начальных численностей $\varepsilon^{-2} \Rightarrow \varepsilon^{-d}$.

**Задача 11 (карикатурная модель стакана; Введенская–Замятин–Малышев)\***. Ось $x$ разбита на ячейки, расстояние между соседними ячейками $\varepsilon$. В ячейках могут находиться частицы "+" и "–". Однако частицы разного знака не могут находиться в одной ячейки, поскольку моментально прореагируют и исчезнут. Каждая частица "+" с интенсивностью $\lambda_+/\varepsilon$ независимо от остальных перепрыгивает в соседнюю справа ячейку (с большей координатой), аналогично, каждая частица "–" с интенсивностью $\lambda_-/\varepsilon$ независимо от стальных перепрыгивает в соседнюю слева ячейку. Введя плотность числа частиц (с учетом знака) $c(t,x)$, считая, что можно корректно определить $c(0,x) = \begin{cases} -2, & x \le 0 \\ 1, & x > 0 \end{cases}$, получите в пределе $\varepsilon \to 0+$ уравнение

$$\frac{\partial c}{\partial t} + \frac{\partial \varphi(c)}{\partial x} = 0, \quad \varphi(c) = \begin{cases} \lambda_+ c, & c > 0 \\ \lambda_- |c|, & c \le 0 \end{cases}.$$

Как понимать решение (начальной) задачи Коши для этого уравнения? Покажите, что обобщенное решение этой задачи Коши будет иметь вид ударной волны. Определите скорость этой ударной волны (формула Римана–Рэнкина–Гюгонио).

**Замечание.** Подробнее об обобщенных решениях законов сохранения, в том числе ударных волнах, написано, например, во второй главе монографии

Введение в математическое моделирование транспортных потоков. М.: МЦНМО, 2013.

К сожалению, из-за негладкости в нуле функции $\varphi(c)$, напрямую нельзя использовать описанные в этой книге наработки. Тем не менее, если немного сгладить функцию $\varphi(c)$, то обобщенное (по Олейник–Кружкову) решение уже будет определено корректно и будет совпадать с формально посчитанной (по формуле Римана–Рэнкина–Гюгонио) ударной волной в первоначальной постановке.

**Задача 12 (игрушечная модель эволюции Д. Кьялво).** В аудитории собрали большое число студентов. Каждый студент, соответствующий определенному биологическому виду, имеет приспособленность (число от нуля до единицы). Динамика (эволюция) состоит в том, что на каждом шаге выбирается наименее приспособленный вид (студент с наименьшей приспособленностью) и еще какие-то два студента, которые выбираются независимо и равновероятно.[14] Для каждого из этих трех студентов случайно (независимо и с равномерным распределением на $[0,1]$) выбираются новые значения приспособленностей. Определим порог приспособленности $c \in [0,1]$. Будем считать, что число студентов стремится к бесконечности, а доля студентов с приспособленностью ниже этого порога стремится к нулю. Тогда $c$ определяется из условия, что после одного шага эволюции среднее число студентов ниже порога приспособленности не изменится. Покажите, что $c = 1/3$. Предположим, что в начальный момент времени все студенты (кроме одного) имеют приспособленности выше пороговой. Лавина начинается тогда когда появляется еще студент(-ы) с приспособленностью ниже пороговой, а заканчивается тогда, когда снова остается всего один студент, имеющий приспособленность выше пороговой. Покажите, что распределение лавин по длительности имеет степенной закон. Здесь имеется в виду, что если понаблюдать за студентами достаточно долго и нарисовать график частота встречаемости лавины / длительность лавины, то с хорошей точностью получится степенной закон с показателем степени $\tau = 3/2$.

**Указание.** Воспользуйтесь степенным законом распределением времени первого возвращения симметричного случайного блуждания на прямой, см., например, *Newman M.E.J.* Power laws, Pareto distributions and Zipf's law // Contemporary Physics. 2005. V. 46. no. 5. P. 323–351. Данная задача немного проясняет один из аспектов того, почему степенные законы так часто возникают при описании зависимости частоты различных катастроф от их масштабы и ряд других зависимостей, встречающихся в "Природе", см. *Бак П.* Как работает природа: Теория самоорганизованной критичности. М.: УРСС, 2013.

**Задача 13 (Majority rule model).** Популяция состоит из $N \gg 1$ агентов, каждый из которых имеет спин (мнение) $\pm 1$. На каждом шаге случайно (независимо и равновероятно) выбираются три агента. Если у всех агентов одинаковый спин ничего не происходит, иначе, большинство (два агента) "уговаривает" заставляет оставшегося агента поменять свой спин. Опишите, что произойдет с агентами по прошествии достаточно большого времени. Постарайтесь получить количественные оценки, в том числе времени выхода в равновесия.

---

[14] В исходной постановке Д. Кьялво студенты были размещены на окружности и "еще два студента" выбирались по принципу соседства с наименее приспособленным. Однако это более сложная для исследования ситуация. Впрочем, можно посчитать: $c \simeq 0.667$, $\tau \simeq 1.07$, но это требует более изощренного определения порога приспособленности.

**Замечание.** Описание более общих моделей можно найти в обзоре *Castellano C., Fortunato S., Loreto V.* Statistical physics of social dynamics // e-print, 2009 arXiv:0710.3256v2. Также как и в предыдущей задаче, модель усложняется, если вводить (пространственно) локальные взаимодействия. В частности, многое зависит от размерности пространства, в котором расположена решетка со спинами. Детали имеются, например, в монографии *Лигетт Т.М.* Марковские процессы с локальным взаимодействием. М.: Мир, 1989.

**Задача 14 (динамика наилучших ответов для игры "Камень-Ножницы-Бумага")\*.** Популяция состоит из $N \gg 1$ агентов. Стратегией каждого агента является выбор камня, ножниц или бумаги. Выигрыш (в единицу времени) агента выбравшего камень есть разность [доля агентов, выбравших ножницы] – [доля агентов, выбравших бумагу], выбравшего ножницы есть разность [доля агентов, выбравших бумагу] – [доля агентов, выбравших камень], выбравшего бумагу есть разность [доля агентов, выбравших камень] – [доля агентов, выбравших ножницы]. При этом с каждым агентом связан свой "пуассоновский будильник"[15]. Все будильники имеют одинаковый параметр. Когда будильник звонит агент выбирает исходя из текущей ситуации камень, ножницы или бумагу, чтобы получить максимальный выигрыш (best response protocol; динамика наилучших ответов). Опишите, что произойдет с агентами по прошествии достаточно большого времени. Постарайтесь получить количественные оценки.

**Замечание.** Для решения этой и следующей задачи рекомендуется посмотреть монографию *Sandholm W.* Population games and Evolutionary dynamics. Economic Learning and Social Evolution. MIT Press; Cambridge, 2010. Интересно отметить, что в этом примере есть чувствительность к выбору "разумной" динамики. Это в некотором смысле нетипично, если речь идет о сходимости к единственному равновесию, как в данной задаче. Скажем, для популярной в популяционной теории игр динамики репликаторов в пределе $N \to \infty$ (этот предел нужен, чтобы можно было смотреть на стохастическую динамику, как на детерминированную) концентрации агентов, использующих разные стратегии, будут колебаться вокруг равновесия $(1/3, 1/3, 1/3)$, и лишь в чезаровском смысле будут сходиться к этому равновесию.

**Задача 15 (модель равновесного распределения транспортных потоков).** Ориентированный граф $\Gamma = (V, E)$ – транспортная сеть города ($V$ – узлы сети (вершины), $E \subset V \times V$ – дуги сети (рёбра графа)). Пусть $W = \{w = (i, j): i, j \in V\}$ – множество пар источник – сток; $p = \{v_1, v_2, ..., v_m\}$ – путь из вершины $v_1$ в $v_m$, если $(v_k, v_{k+1}) \in E$, $k = 1,...,m-1$; $P_w$ – множество путей, отвечающих корреспонденции $w \in W$; $P = \bigcup_{w \in W} P_w$ – совокупность всех путей в сети $\Gamma$; $x_p$ – величина потока по пути $p$, $x = \{x_p : p \in P\}$; $y_e$ – величина потока по дуге $e$: $y_e = \sum_{p \in P} \delta_{ep} x_p$, где $\delta_{ep} = \begin{cases} 1, & e \in p \\ 0, & e \notin p \end{cases}$ ($y = \Theta x$, $\Theta = \{\delta_{ep}\}_{e \in E, p \in P}$); $\tau_e(y_e)$ – удельные затраты на проезд по дуге $e$ (гладкие неубывающие функции); $G_p(x)$ $G_p(x) = \sum_{e \in E} \tau_e(y_e) \delta_{ep}$ – удельные затраты на проезд по пути $p$, $G(x) = \{G_p(x) : p \in P\}$ ($G(x) = \Theta^T \tau(y)$). Пусть известна матрица корреспонденций $\{d_w\}_{w \in W}$. Тогда вектор $x$, характеризующий распределение потоков, должен лежать в допустимом множестве: $X = \left\{ x \geq 0 : \sum_{p \in P_w} x_p = d_w, w \in W \right\}$. Рассмотрим игру, в которой каждому элементу $w \in W$ соответствует свой, достаточно большой, набор однотипных агентов ("сидящих на корреспонденции $w$"). Множеством чистых стратегий каждого такого агента является $P_w$, а выигрыш (потери со знаком минус) от использования стратегии $p \in P_w$ определяются формулой $-G_p(x)$.

---

[15] Пуассоновский будильник звонит через случайные (одинаково распределенные) промежутки времени, каждый промежуток – независимая показательная случайная величина. Параметр этой случайной величины называется параметром будильника. Фактически мы уже сталкивались с пуассоновским будильником, например, в задаче "парадокс Эренфестов" и в ее распределенном варианте.

Пусть $d_w := d_w M$ ($M \gg 1$), $x := x/M$, $y := y/M$, $\tau_e(y_e) := \tau_e(y_e/M)$.

**а)\*** Покажите, что при $M \to \infty$ динамика наилучших ответов из предыдущей задачи "приводит" к равновесию $x^*$ Нэша–Вардропа, которое для данной игры загрузок (а, стало быть, потенциальной игры) сводится к задаче выпуклой оптимизации (потенциала):

$$\Psi(y(x)) = \sum_{e \in E} \int_0^{y_e(x)} \tau_e(z) dz \to \min_{x \in X}.$$

**б) (энтропийная регуляризация)\*\*** Легко проверить, что если $\tau'_e(\cdot) > 0$, $e \in E$, то равновесное распределение потоков по рёбрам $y^*$ (см. п. а)) единственно. Но это еще не гарантирует единственности $x^*$ (напомним связь: $y(x) = \Theta x$). Когда равновесных распределений потоков по путям может быть много, то нужно понять, а какое из этих равновесий все-таки реализуется на практике (если реализуется, ведь может быть, например, ни к чему не сходящееся блуждание по равновесному множеству). Предположим, что каждый агент независимо принимает решения на основе текущей зашумленной информации $\arg\max_{p \in P_w} \{-G_p(x) + \xi_p\}$, где $\xi_p$ – независимые одинаково распределённые случайные величины с нулевым математическим ожиданием и ограниченной дисперсией. Будем считать, что $\xi_p$ имеют распределением Гумбеля: $P(\xi_p < \xi) = \exp\{-e^{-\xi/\sigma - E}\}$, $E \approx 0.5772$ – константа Эйлера, а $D\xi_p = \sigma^2 \pi^2/6$. Покажите, что при $M \to \infty$ такая динамика (её называют имитационная логит динамика) приводит к, так называемому, стохастическому равновесию, которое определяется, как единственное решение сильно выпуклой (в 1-норме) задачи оптимизации:

$$\Psi(y(x)) + \sigma \sum_{w \in W} \sum_{p \in P_w} x_p \ln(x_p/d_w) \to \min_{x \in X}.$$

Если $\sigma \to 0+$, то решение этой задачи сходится к

$$x^* = \arg\min_{x \in X : \Theta x = y^*} \sum_{w \in W} \sum_{p \in P_w} x_p \ln(x_p/d_w),$$

где $y^* = \arg\min_{y = \Theta x, \, x \in X} \Psi(y)$.

**Замечание.** Распределение потоков по путям $x = \{x_p\} \in X$ называется равновесием (Нэша–Вардропа) в популяционной игре $\langle \{x_p\} \in X, \{G_p(x)\} \rangle$, если из $x_p > 0$ ($p \in P_w$) следует $G_p(x) = \min_{q \in P_w} G_q(x)$. Или, что то же самое:

для любых $w \in W$, $p \in P_w$ выполняется $x_p \cdot \left( G_p(x) - \min_{q \in P_w} G_q(x) \right) = 0$.

В этой задаче имеется ряд моментов, которые связаны с возможностью осуществления предельного перехода по числу агентов стремящихся к бесконечности, и требуют некоторых дополнительных оговорок. Детали имеются, например, в уже упоминавшейся книге W. Sandholma (в особенности, стоит обратить внимание на теорему 11.5.12). Эта задача может быть обобщена на многостадийные транспортные модели, которые сейчас популярны на практике:

*Гасников А.В., Дорн Ю.В., Нестеров Ю.Е., Шпирко С.В.* О трехстадийной версии модели стационарной динамики транспортных потоков // Математическое моделирование. 2014. Т. 26. № 6. С. 34–70. arXiv:1405.7630

*Гасников А.В.* Об эффективной вычислимости конкурентных равновесий в транспортно-экономических моделях // Математическое моделирование. 2015. Т. 27. № 12. С. 121–136. arXiv:1410.3123

*Бабичева Т.С., Гасников А.В., Лагуновская А.А., Мендель М.А.* Двухстадийная модель равновесного распределения транспортных потоков // Труды МФТИ. 2015. Т. 7. № 3. С. 31–41. https://mipt.ru/upload/medialibrary/971/31-41.pdf

*Гасников А.В., Лагуновская А.А., Морозова Л.Э.* О связи имитационной логит динамики в популяционной теории игр и метода зеркального спуска в онлайн оптимизации на примере задачи выбора кратчайшего маршрута // Труды МФТИ. 2015. Т. 7. № 4. С. 104–113. arXiv:1511.02398

*Гасников А.В., Гасникова Е.В., Мациевский С.В., Усик И.В.* О связи моделей дискретного выбора с разномасштабными по времени популяционными играми загрузок // Труды МФТИ. 2015. Т. 7. № 4. С. 129–142. arXiv:1511.02390

*Баймурзина Д.Р., Гасников А.В., Гасникова Е.В.* Теория макросистем с точки зрения стохастической химической кинетики // Труды МФТИ. 2015. Т. 7. № 4. С. 95–103. https://mipt.ru/upload/medialibrary/ae2/95-103.pdf

*Гасников А.В., Гасникова Е.В., Двуреченский П.Е., Ершов Е.И., Лагуновская А.А.* Поиск стохастических равновесий в транспортных моделях равновесного распределения потоков // Труды МФТИ. 2015. Т. 7. № 4. С. 114–128. arXiv:1505.07492

*Гасников А.В., Гасникова Е.В., Мендель М.А., Чепурченко К.В.* Эволюционные выводы энтропийной модели расчета матрицы корреспонденций // Математическое моделирование. 2016. Т. 28. № 4. С. 111–124. arXiv:1508.01077

Описанная в этой задаче зашумленная динамика наилучших ответов играет важную роль при описании поведения ограничено рациональных агентов:

*Andersen S.P., de Palma A., Thisse J.-F.* Discrete choice theory of product differentiation. MIT Press; Cambridge, 1992;
*Sandholm W.* Population games and Evolutionary dynamics. Economic Learning and Social Evolution. MIT Press; Cambridge, 2010.

**Задача 16 (эволюция РНК по В.Г. Редько)\*.** Для простоты рассуждений будем считать, что молекула РНК $S$ представляет собой цепочку битов длины $N \gg 1$. В системе имеется $n$ ($N \ll n \ll 2^N$) молекул РНК. Для определенности, в начальный момент все молекулы приготовлены независимо и случайно (из равномерного распределения). Будем считать, что самая приспособленная молекула $S_0$ – это молекула, состоящая из одних нулевых битов. Ее приспособленность $f(S_0) = 1$. Приспособленность молекулы $S$ – $f(S) = \exp(-\rho(S, S_0))$, где $\rho$ – расстояние Хэмминга. Динамика системы молекул в дискретном времени заключается в следующем. Отбираем $n$ особей в новую популяцию, делая $n$ независимых разыгрываний из дискретного распределения с вероятностями соответствующих исходов $\sim \{f(S_k)\}_{k=1}^{n}$ (естественный отбор). Далее в получившейся популяции подвергаем независимо каждый бит каждой молекулы случайной мутации. Вероятность мутации $P \sim 1/N$. Покажите, что после $T = \mathrm{O}(N)$ шагов произойдет "фиксация" наиболее приспособленного вида. Другими словами с большой вероятностью большая часть молекул будет иметь почти все свои биты нулевыми. Постарайтесь количественно описать последнее предложение.

**Замечание.** По поводу постановки задачи и ее окрестностей см. *Редько В.Г.* Эволюция, нейронные сети, интеллект: Модели и концепции эволюционной кибернетики. М.: УРСС, 2013. Для оценки скорости сходимости (mixing time) в этой и других задачах рекомендуется использовать монографию *Levin D.A., Peres Y., Wilmer E.L.* Markov chain and mixing times. AMS, 2009. http://pages.uoregon.edu/dlevin/MARKOV/markovmixing.pdf Если функция приспособленности многоэкстремальная, а требуется найти глобальный максимум, то описанную в задаче динамику можно понимать как один из вариантов генетического алгоритма, отыскивающего глобальный максимум. При этом в генетических алгоритмах стараются выбирать $n$ на границе условия $N \ll n$, часто даже $n \simeq N$. Хотя генетические алгоритмы достаточно популярны в задачах глобальной оптимизации (см., например, *Zhigljavsky A., Zilinskas A.* Stochastic global optimization. Springer Optimization and Its Applications, 2008), однако об их сходимости (точнее об оценках скорости сходимости) известно не так много. Впрочем, некоторые результаты все-таки есть, см., например, статью *Cerf R.* Asymptotic convergence of genetic algorithms // Adv. Appl. Prob. 1998. V. 30. no. 2. P. 521–550.

**Задача 17 (эволюция РНК с рекомбинациями, но без отбора). а)\*** Для простоты рассуждений будем считать, что молекула РНК $S$ представляет собой цепочку битов длины $N$. В системе имеется $n \gg 2^N$ молекул РНК. Для определенности, в начальный момент все молекулы приготовлены независимо и случайно (из равномерного распределения). С интенсивностью $\lambda$ каждый 0 бит каждой молекулы независимо от остальных превращается в 1 бит, а с интенсивностью $\mu$ – наоборот.[16] Считаем также, что каждой паре молекул $x$ и $y$ соответствует $2^N$ однотипных реакций (рекомбинаций). Реакция заключается в том, что какой-то фрагмент $x_I$ (цепочка битов, не обязательно идущих подряд) одной молекулы копируется и заменяет соответствующий участок $y_I$ в другой молекуле. Интенсивность таких реакций мы считаем известной $0 \le \varphi(x_I, y_I) \le 1$, причем $\varphi(x_I, y_I) = \varphi(y_I, x_I)$. Обозначим через $\mu_k(t)$ – долю молекул с числом единичных битов равным $k$. Опишите поведение системы молекул в терминах $\mu_k(t)$ на больших временах. Отметим, что это поведение не зависит от вида функции $\varphi(x_I, y_I)$. Рассмотрите случай, когда $\lambda = 1$, $\mu = 2$. Постарайтесь получить количественные оценки.

**б)\*\*** Попробуйте предложить такое обобщение модели из п. а), в котором бы учитывался естественный отбор (или такое обобщение модели из предыдущей задачи, в котором учитывалась бы рекомбинация), и которое позволяет аналитически исследовать поведение на больших временах.

**Замечание.** Для лучшего понимания постановки этой задачи рекомендуется посмотреть цикл недавних работ S.A. Pirogov, A.N. Rybko и др. на arxiv.org. Для решения данной задачи полезно выписать с помощью теоремы Т. Куртца (впрочем, эта теорема полезна практически во всех "эволюционных" задачах) систему $2^N$ нелинейных дифференциальных уравнения на $\mu^x(t)$, где $\mu^x(t)$ – доля молекул с набором битов $x$. Эта система получает при предельном переходе $n \to \infty$, в предположении, что в начальный момент существуют соответствующие предельные пропорции $\mu^x(0)$. Детали см. в *Ethier N.S., Kurtz T.G.* Markov processes. Wiley Series in Probability and Mathematical Statistics: Probability and Mathematical Statistics. John Wiley & Sons Inc., New York, 2005.

**Задача 18 (Пуассоновская гипотеза)\*.** Имеется $N$ заявок и $M$ серверов с процедурой обслуживания FIFO. После обслуживания каждая заявка с равной вероятностью отправляется на один из серверов. Времена обслуживания заявок – независимые одинаково распределенные с.в. с ф.р. $F(x)$, $\mu = \int_0^\infty x\,dF(x)$. Предположим, что $N \to \infty$, $M/N \to \lambda$. Исследуйте поведение системы на больших временах.

---

[16] Можно сказать так, что с каждой молекулой может происходить $N$ различных типов реакций, в результате реакции меняется один из битов (соответствующей данной реакции). Каждой молекуле и каждому типу реакции соответствует свой пуассоновский будильник. Когда будильник звонит, происходит соответствующая реакция.

**Замечание.** Данная задача является частным случаем более общего класса задач, в которых изучаются сети (массового обслуживания) Джексона (и их обобщения) при термодинамическом предельном переходе, см., например, цикл работ A. Rybko, S. Shlosman на http://arxiv.org/. Отметим в этой связи интересные исследования фазового перехода в таких сетях (см. приложение А.А. Замятина, В.А. Малышева в книге Введение в математическое моделирование транспортных потоков. М: МЦНМО, 2013). Однако у приведенной задачи есть особенность – время обслуживание, вообще говоря, не предполагается экспоненциальным. Вообще эта задача очень показательна во многих отношениях (см. статью А.Н. Рыбко в 4-м номере "Глобус").

**Задача 19 (теорема Гордона–Ньюэлла и PageRank).** а) Имеется $N \gg 1$ пользователей, которые случайно (независимо) блуждают в непрерывном времени по ориентированному графу с эргодической инфинитезимальной матрицей $\Lambda$. Назовем вектор $p$ (из единичного симплекса) PageRank, если $\Lambda p = 0$. Обозначим через $n_i(t)$ – число пользователей на $i$-й странице в момент времени $t \geq 0$. Покажите, что $n(t)$ асимптотически имеет мультиномиальное распределение с вектором параметров PageRank $p$, т.е.

$$\lim_{t \to \infty} P(n(t) = n) \sim \prod_i \frac{(p_i)^{n_i}}{n_i!}.$$

Следовательно (неравенство Хефдинга в гильбертовом пространстве),

$$\lim_{t \to \infty} P\left( \left\| \frac{n(t)}{N} - p \right\|_2 \geq \frac{2\sqrt{2} + 4\sqrt{\ln(\sigma^{-1})}}{\sqrt{N}} \right) \leq \sigma.$$

**б)** Получите тот же результат, что в п. а) рассмотрев соответствующую систему унарных химических реакций. Переход одного из пользователей из вершины $i$ в вершину $j$ – означает превращение одной молекулы вещества $i$ в одну молекулу вещества $j$, $n_i(t)$ – число молекул $i$-го типа в момент времени $n_i(t)$. Каждое ребро графа соответствует определенной реакции (превращению). Интенсивность реакций определяется матрицей $\Lambda$ и числом молекул, вступающих в реакцию (закон действующих масс). Покажите, что условие $\Lambda p = 0$ – в точности соответствует условию унитарности в стохастической химической кинетике (обобщению условия детального баланса, предложенного в 2000 Батищевой–Веденяпиным и независимо Малышевым–Пироговым–Рыбко в 2004).

**Замечание.** По п а) полезно посмотреть монографию *Serfozo R.* Introduction to stochastic networks. Springer, 1999, а по п. б) статью *Гасников А.В., Гасникова Е.В.* Об энтропийно-подобных функционалах … // Математические заметки. 2013. Т. 94. № 6. С. 816–824 и цитированную там литературу. С примером более общих реакций можно познакомиться, например, по задаче "Кинетика социального неравенства", в которой реакции бинарные. Отметим, что если воспользоваться теоремой Санова о больших уклонениях для мультиномиального распределения, то получим

$$\frac{N!}{n_1!\cdot\ldots\cdot n_m!}p_1^{n_1}\cdot\ldots\cdot p_m^{n_m} = \exp\left(-N\sum_{i=1}^{m}v_i\ln(v_i/p_i)+R\right),$$

где $v_i = n_i/N$, $|R| \leq \frac{m}{2}(\ln N + 1)$. Однако последующее применение неравенства Пинскера не дает нам равномерной по $m$ оценки в 1-норме. Как и ожидалось, выписанная в задаче оценка в 2-норме (правая часть неравенства под вероятностью) и так полученная оценка в 1-норме будут отличаться по порядку приблизительно в $\sqrt{m}$ раз, что соответствует типичному соотношению между 1 и 2 нормами. Слово "типично" здесь отвечает, грубо говоря, за ситуацию, когда компоненты вектора одного порядка. Для многих приложений, где возникают предельные конфигурации (кривые), описывающиеся вектором с огромным числом компонент, имеет место быстрый закон убывания этих компонент. Например, для некоторого обобщения (см. *Райгородский А.М.* Модели интернета. Долгопрудный: Изд. Дом "Интеллект", 2013) модели Барабаши–Альберт случайного роста графа интернета, вектор PageRank с хорошей точностью и с высокой вероятностью имеет компоненты, убывающие по степенному закону. Если в этих приложениях возникает мультиномиальное распределение с таким вектором, то при изучении концентрации этого распределения предпочтительнее становится подход с 1-нормой (см., например, задачу 4 "Кинетика социального неравенства").

К п. б) можно привести следующее обобщение. Предположим, что некоторая макросистема может находиться в различных состояниях, характеризуемых вектором $n$ с неотрицательными целочисленными компонентами. Будем считать, что в системе происходят случайные превращения (химические реакции).

Пусть $n \to n - \alpha + \beta$, $(\alpha, \beta) \in J$ – все возможные типы реакций, где $\alpha$ и $\beta$ – вектора с неотрицательными целочисленными компонентами. Введем интенсивность реакции:

$$\lambda_{(\alpha,\beta)}(n) = \lambda_{(\alpha,\beta)}(n \to n-\alpha+\beta) = N^{1-\sum_i \alpha_i} K_\beta^\alpha \prod_{i:\alpha_i>0} n_i \cdot \ldots \cdot (n_i - \alpha_i + 1),$$

где $K_\beta^\alpha \geq 0$ – константа реакции; при этом $\sum_{i=1}^{m} n_i(t) \equiv N \gg 1$. Другими словами, $\lambda_{(\alpha,\beta)}(n)$ – вероятность осуществления в единицу времени перехода $n \to n - \alpha + \beta$. Здесь не предполагается, что число состояний $m = \dim n$ и число реакций $|J|$ не зависят от числа агентов $N$. Тем не менее, если ничего не известно про равновесную конфигурацию $c^*$ (типа быстрого убывания компонент этого вектора), то дополнительно предполагается, что $m \ll N$ – это нужно для обоснования возможности применения формулы Стирлинга при получении вариационного принципа (максимума энтропии) для описания равновесия макросистемы $c^*$ (в концентрационной форме). Однако часто априорно можно предполагать (апостериорно проверив), что компоненты вектора $c^*$ убывают быстро, тогда это условие можно отбросить. Так, например, обстоит дело с "Кинетикой социального неравенства".

Возникающий марковский процесс считается неразложимым. Далее приводится **теорема** (во многом установленная в 1999 – 2005 гг. в работах Я.Г. Батищевой, В.В. Веденяпина, В.А. Малышева, С.А. Пирогова, А.Н. Рыбко).

**а)** $\langle \mu, n(t) \rangle \equiv \langle \mu, n(0) \rangle$ (inv) $\Leftrightarrow$ *вектор $\mu$ ортогонален каждому вектору семейства* $\{\alpha - \beta\}_{(\alpha,\beta) \in J}$. *Здесь* $\langle \cdot, \cdot \rangle$ – *обычно евклидово скалярное произведение.*

**б)** (Т. Куртц) *Если существует* $\lim_{N \to \infty} n(0)/N = c(0)$, $K_\beta^\alpha := K_\beta^\alpha(n/N)$, $m$ *и* $|J|$ *не зависит от* $N$, *то для любого* $t > 0$ *с вероятностью 1 существует* $\lim_{N \to \infty} n(t)/N = c(t)$, *где* $c(t)$ – *не случайная вектор-функция, удовлетворяющая СОДУ Гульдберга–Вааге:*

$$\frac{dc_i}{dt} = \sum_{(\alpha,\beta) \in J} (\beta_i - \alpha_i) K_\beta^\alpha(c) \prod_j c_j^{\alpha_j}. \tag{ГВ}$$

*Гиперплоскость (inv) (с очевидной заменой $n \Rightarrow c$) инвариантна относительно этой динамики.*[17] *Более того, случайный процесс $n(t)/N$ слабо сходится при $N \to \infty$ к $c(t)$ на любом конечном отрезке времени.*

**в)** *Пусть выполняется условие унитарности (очевидно, что $\xi$, удовлетворяющий условию (U), – неподвижная точка в (ГВ))*

$$\exists\, \xi > 0 : \forall\, \beta \to \sum_{\alpha:(\alpha,\beta) \in J} K_\beta^\alpha \prod_j (\xi_j)^{\alpha_j} = \sum_{\alpha:(\alpha,\beta) \in J} K_\alpha^\beta \prod_j (\xi_j)^{\beta_j}. \tag{U}$$

*Тогда неотрицательный ортант $\mathbb{R}_+^m$ расслаивается гиперплоскостями (inv), так что в каждой гиперплоскости (inv) уравнение (U) (положительно) разрешимо притом единственным образом. Стало быть, существует, притом единственная, неподвижная*

---

[17] В каком-то смысле жизнь нелинейной динамической системы определяется линейными законами сохранения, унаследованными ею при скейлинге (каноническом). Этот тезис имеет, по-видимому, более широкое применение (см. работы В.В. Веденяпина).

*точка $c^* \in (inv)$ у системы (ГВ), являющаяся при этом глобальным аттрактором. Система (ГВ) имеет функцию Ляпунова $KL(c,\xi) = \sum_{i=1}^{m} c_i \ln(c_i/\xi_i)$.*

*Стационарное (инвариантное) распределение описанного марковского процесса имеет носителем множество (inv) и (с точностью до нормирующего множителя) имеет вид*

$$\frac{N!}{n_1! \cdot \ldots \cdot n_m!} (\xi_1)^{n_1} \cdot \ldots \cdot (\xi_m)^{n_m} \sim \exp(-N \cdot KL(c,\xi)),$$

*где $\xi$ – произвольное решение (U), не важно какое именно (от этого, конечно, будет зависеть нормирующий множитель, но это ни на чем не сказывается). При этом условие унитарности (U) является обобщением условия детального равновесия[18] (баланса)[19]*

$$\exists\ \xi > 0:\ \forall\ (\alpha,\beta) \in J \to K_\beta^\alpha \prod_j (\xi_j)^{\alpha_j} = K_\alpha^\beta \prod_j (\xi_j)^{\beta_j},$$

*принимающего такой вид для мультиномиальной стационарной меры.*

*Используя неравенство Чигера (для оценки mixing time) и неравенство Хефдинга в гильбертовом пространстве, отсюда можно получить (зависимость $a(m,c(0))$ во многих приложениях может быть равномерно ограничена)*

$$\exists\ a = a(m,c(0)):\ \forall\ \sigma > 0, t \geq a \ln N \to P\left( \left\| \frac{n(t)}{N} - c^* \right\|_2 \geq \frac{2\sqrt{2} + 4\sqrt{\ln(\sigma^{-1})}}{\sqrt{N}} \right) \leq \sigma,$$

---

[18] В терминах условия задачи условие унитарности просто означает, что в равновесии для любой вершины имеет место баланс числа пользователей входящих в единицу времени в эту вершину с числом пользователей, выходящих в единицу времени из этой вершины. В то время как условие детального равновесия означает, что в равновесии для любой пары вершин число пользователей, переходящих в единицу времени из одной вершину в другую равно числу пользователей, переходящих в обратном направлении. Понятно, что второе условие является частным случаем первого.

[19] Много интересных примеров макросистем, в которых $K_\beta^\alpha := K_\beta^\alpha(n/N)$, и имеет место детальный баланс, собрано в книге *Вайдлих В.* Социодинамика: системный подход к математическому моделированию в социальных науках. М.: УРСС, 2010. Выше (задача 15) нами был рассмотрен пример на эту тему из книги W. Sandholma. Как правило, такие постановки приводят к функциям Ляпунова–Санова вида $\Psi(c) + \sigma KL(c,\xi)$, $\sigma \geq 0$. Причем во всех этих приложениях условие детального баланса выполняется точно. Отметим, что с некоторыми дополнительными оговорками можно допускать и приближенное выполнение условий детального баланса (с точностью $O(N^{-1})$).

*где (принцип максимума энтропии Больцмана–Джейнса)*

$$c^* = \arg\max_{c \in (\text{inv})} \left( -\sum_i c_i \ln(c_i/\xi_i) \right) = \arg\min_{c \in (\text{inv})} KL(c, \xi),$$

*а $\xi$ – произвольное решение (U), причем $c^*$ определяется единственным образом, т.е. не зависит от выбора $\xi$. Геометрически $c^*$ – это KL-проекция произвольного $\xi$, удовлетворяющего (U), на гиперплоскость (inv), соответствующую начальным данным $c(0)$. Независимость этой проекции от выбора $\xi$ из (U) просто означает, что кривая (U) проходит KL-перпендикулярно через множество (inv).*

**г)** (Е.В. Гасникова) *Верно и обратное утверждение, то есть условие (U) не только достаточное для того чтобы равновесие находилось из приведенной выше задачи энтропийно-линейного программирования, но и, с небольшой оговоркой (для почти всех $c(0)$), необходимое. Также верно и более общее утверждение, связывающее понимание энтропии в смысле Больцмана (функция Ляпунова прошкалированной кинетической динамики) и Санова (функционал действия в неравенствах больших уклонений для инвариантной меры): если стационарная мера асимптотически представима в виде $\sim \exp(-N \cdot V(c))$, то $V(c)$ – функция Ляпунова (ГВ).*

К п. г) можно сделать следующее пояснение. Обозначим через $h(c)$ вектор-функцию, стоящую в правой части СОДУ (ГВ). Тогда (см., например, книгу К.В. Гардинера) при $N \gg 1$ по теореме Т. Куртца $n(t)/N$ будет $\mathrm{O}\left(\log N/\sqrt{N}\right)$-близко (в этом месте требуются оговорки, чтобы близость была равномерна по времени) к $x_t = x(t)$ – решению стохастической системы дифференциальных уравнений (с начальным условием $x_0 = c(0)$)

$$dx_t = h(x_t)dt + \sqrt{\frac{g(x_t)}{N}}dW_t,$$

где функция $g(x_t) > 0$ рассчитывается по набору реакций и константам реакций (которые могут быть не постоянны и зависеть от концентраций). Инвариантная (стационарная) мера $m(x) = \lim_{t \to \infty} p(t, x)$ этого однородного марковского процесса удовлетворяет уравнению (см. также задачу 21)

$$\frac{1}{2N}\nabla^2\big(g(x)m(x)\big) - \operatorname{div}\big(h(x)m(x)\big) = 0,$$

поскольку плотность распределения $p(t,x)$ процесса $x_t$ подчиняется уравнению Колмогорова–Фоккера–Планка

$$\frac{\partial p(t,x)}{\partial t} = -\operatorname{div}\big(h(x)p(t,x)\big) + \frac{1}{2N}\nabla^2\big(g(x)p(t,x)\big).$$

Здесь использовалось обозначение: $\nabla^2 f(x) = \sum_{i,j} \partial^2 f(x)/\partial x_i \partial x_j$. Если известно, что

$$m(x) \simeq \operatorname{const} \cdot \exp\big(-N \cdot V(x)\big),$$

то из уравнения на $m(x)$ имеем

$$N\langle h, \nabla V\rangle - \operatorname{div} h - \frac{1}{2}\langle \nabla g, 1\rangle\langle \nabla V, 1\rangle + \frac{1}{2N}V\nabla^2 g - \frac{1}{2}g\nabla^2 V + \frac{N}{2}g\langle \nabla V, 1\rangle^2 \simeq 0,$$

следовательно

$$\langle h, \nabla V\rangle \simeq -\frac{1}{2}g\langle \nabla V, 1\rangle^2 + \mathrm{O}\!\left(\frac{1}{N}\right)\underset{N\to\infty}{=} -\frac{1}{2}g\langle \nabla V, 1\rangle^2 \le 0.$$

Эта выкладка поясняет (но не доказывает, для доказательства требуются более аккуратные рассуждения), почему функция $V(c)$ может быть функцией Ляпунова системы (ГВ) $dc/dt = h(c)$. Более того, модель стохастической химической кинетики здесь может быть заменена и более общими шкалирующимися марковскими моделями.

**Задача 20 (теорема Гордона–Ньюэлла; Л.Г. Афанасьева)\*.** Рассматривается транспортная сеть, в которой между $N$ станциями курсируют $M$ такси. Клиенты пребывают в $i$-й узел в соответствии с пуассоновским потоком с параметром $\lambda_i > 0$ ($i = 1, ..., N$). Если в момент прибытия в $i$-й узел там есть такси, клиент забирает его и с вероятностью $p_{ij} \ge 0$ направляется в $j$-й узел, по прибытии в который покидает сеть. Такси остается ждать в узле прибытия нового клиента. Времена перемещений из узла в узел – независимые случайные величины, имеющие показательное распределение с параметром $v_{ij} > 0$ для пары узлов $(i, j)$. Если в момент прихода клиента в узел там нет такси, клиент сразу покидает узел. Считая $p_{ij} = N^{-1}$, $\lambda_i = \lambda$, $v_{ij} = v$, покажите, что вероятность того, что клиент, поступивший в узел (в установившемся (стационарном) режиме работы сети), получит отказ, равна

$$p_{\text{отказа}}(N, M) = \sum_{k=0}^{M} \frac{C_{N-2+k}^{k} \rho^{M-k}}{(M-k)!} \bigg/ \sum_{k=0}^{M} \frac{C_{N-1+k}^{k} \rho^{M-k}}{(M-k)!}, \ \rho = N\lambda/\nu.$$

Методом перевала покажите справедливость следующей асимптотики при $N \to \infty$:

$$p_{\text{отказа}}(N, rN) = 1 - \frac{2r}{\lambda/\nu + r + 1 + \sqrt{(\lambda/\nu + r + 1)^2 - 4\lambda r/\nu}} + \text{O}\left(\frac{1}{N}\right).$$

**Замечание.** Метод перевала – очень полезный инструмент исследования асимптотик интегралов по параметру. Его подробное изложение имеется, например, здесь: *Федорюк М.В.* Метод перевала. М.: УРСС, 2010.

**Задача 21 (Mean field games; Lasry–Lions)**\*\*. Пусть динамика агента (игрока) задается стохастическим дифференциальным уравнением

$$dX_t = -\alpha_t dt + \sigma dW_t, \ X_0 = x_0$$

где $W_t$ – винеровский процесс, а $\alpha_t$ – стратегия (управление), $\alpha_t, X_t \in \mathbb{R}^d$, Функционал потерь задается формулой

$$J(x_0, \alpha) = \lim_{T \to \infty} \frac{1}{T} E\left[\int_0^T \left(L(\alpha_t) + F(X_t)\right) dt\right],$$

где $L(\alpha)/|\alpha| \to \infty$ при $|\alpha| \to \infty$. Считаем, что функция $F(x)$ 1-периодическая по каждой компоненте вектора $x$. Таким образом, можно ограничится рассмотрением стохастической динамики на $\mathrm{T}^d$ – единичном торе в $\mathbb{R}^d$. Введем гамильтониан $H(p) = \sup_{\alpha \in \mathbb{R}^d} \{\langle p, \alpha \rangle - L(\alpha)\}$ и $\nu = \sigma^2/2$.

**а) (Fleming–Soner; Bardi–Capuzzo-Dolcetta)** Покажите, что если уравнение Гамильтона–Якоби–Беллмана имеет решение (пара $\lambda$, $V(\cdot)$)

$$-\nu\Delta V + H(\nabla V) + \lambda = F(x),$$

тогда оптимальное значение функционала $\inf_{\alpha} J(x_0, \alpha) = J(x_0, \hat{\alpha}) = \lambda$ и оптимальное управление (в форме синтеза) $\hat{\alpha}(x) = \nabla H(\nabla V(x))$. Кроме того, марковский (диффузионный) процесс $dX_t = -\hat{\alpha}_t(X_t) dt + \sigma dW_t$ имеет инвариантную меру $m(x) dx$, являющуюся решением уравнения Колмогорова–Фоккера–Планка

$$\nu\Delta m + \text{div}(\nabla H(\nabla V) m) = 0, \ \int_{\mathrm{T}^d = [0,1]^d} m(x) dx = 1, \ m(x) > 0.$$

Ввиду эргодичности $X_t$, $m(x)$ – асимптотическое распределение вероятностей положения агента, ведущего себя оптимальным образом.

**б) (A. Friedman; Bensoussan–Frehse)** Предположим теперь, что имеется $N>1$ игроков, и каждого своя динамика

$$dX_t^i = -\alpha_t^i dt + \sigma^i dW_t, \ X_0^i = x_0^i,$$

своя функция потерь

$$J^i(\{x_0\},\{\alpha\}) = \varlimsup_{T\to\infty}\frac{1}{T} E\left[\int_0^T \left(L^i(\alpha_t^i) + F^i(X_t^1,...,X_t^N)\right)dt\right].$$

Аналогично видоизменяются определения $H^i(p)$ и $v^i$. Равновесие Нэша $\{\bar{\alpha}\}$ в такой игре определяется из условия (для любого $i=1,...,N$)

$$J^i(\{x_0\},\{\bar{\alpha}\}) = \min_{\alpha^i} J^i\left(\{x_0\},\bar{\alpha}^1,...,\bar{\alpha}^{i-1},\alpha^i,\bar{\alpha}^{i+1},...,\bar{\alpha}^N\right).$$

Можно показать, что система уравнений ($i=1,...,N$)

$$-v^i \Delta V_i + H^i(\nabla V_i) + \lambda_i = f^i(x; m_1,...,m_N),$$

$$v^i \Delta m_i + \operatorname{div}\left(\nabla H^i(\nabla V_i) m_i\right) = 0,$$

$$\int_{\mathrm{T}^d} m_i(x)dx = 1, \ m_i(x) > 0, \ \int_{\mathrm{T}^d} V_i(x)dx = 0,$$

где

$$f^i(x; m_1,...,m_N) = \int_{T^{d(N-1)}} F^i\left(x^1,...,x^{i-1},x,x^{i+1},...,x^N\right) \prod_{j\neq i} dm_j(x^j),$$

имеет решение $\lambda_i$, $V_i(\cdot)$, $m_i(\cdot)$, $i=1,...,N$. Покажите, что для любого такого решения: $\hat{\alpha}^i(x) = \nabla H^i(\nabla V_i(x))$ – равновесие Нэша и $\lambda_i = J^i(\{x_0\},\{\hat{\alpha}\})$ – соответствующее значение игры.

**в)** Пусть для всех игроков $v^i = v$, $H^i = H$, а $F^i$ зависит только от $X^i$ и эмпирической плотности распределения других игроков, т.е.

$$F^i(x^1,...,x^N) = \Phi\left[\frac{1}{N-1}\sum_{j\neq i}\delta_{x_j}\right](x^i),$$

где

$$\Phi: \{\text{распределения вероятностей на } \mathrm{T}^d\} \to \{\text{липшицевы функции на } \mathrm{T}^d\}$$

– непрерывный интегральный функционал

$$\Phi[m](x) = G\left(x, \int_{T^d} k(x,y) dm(y)\right).$$

Можно показать, что последовательность $\lambda_i^{(N)}$, $V_i^{(N)}(\cdot)$, $m_i^{(N)}(\cdot)$ ($N \to \infty$) относительно компактна в $R \times C^2(T^d) \times W^{1,p}(T^d)$ (для любого $1 \le p < \infty$). Следовательно, можно выбрать сходящуюся подпоследовательность. Зафиксируем $i$. Покажите, что при $N \to \infty$ предел любой такой сходящейся подпоследовательности является решением системы ($x \in T^d$)

$$-\nu\Delta V + H(\nabla V) + \lambda = \Phi[m](x),$$

$$\nu\Delta m + \text{div}(\nabla H(\nabla V) m) = 0,$$

$$\int_{T^d} m(x) dx = 1, \; m(x) > 0, \; \int_{T^d} V(x) dx = 0.$$

Можно показать, что решение этой системы существует. Если $\Phi[m](x) = G(m(x))$, $G$ – возрастающая функция, а $H$ – сильно выпуклая функция, то можно гарантировать и единственность.

**Замечание.** Игры среднего поля – относительно новая (возраст около 10 лет) и бурно развивающаяся сейчас область исследований. В России, насколько нам известно, в этом направлению работают небольшие группы в Екатеринбурге (Ю.В. Авербух) и в Москве–Лондоне (В.Н. Колокольцов). Очень бурно это направления исследуется на Западе французскими и американскими математиками: P.L. Lions, J.M. Lasry, O. Gueant, M. Bardi, A. Bensoussan, J. Frehse, D. Lacker и др.

**Задача 22 (круговая модель М. Каца).** На окружности отмечено $n$ равноотстоящих точек, $m$ из которых составляют множество $Q$. В каждой точке помещен белый (б) или черный (ч) шар. В единицу времени каждый шар переходит против часовой стрелки в соседнюю точку, причем, если шар выходит из точки, принадлежащей $Q$, он меняет свой цвет. Пусть $\mu = m/n < 1/2$. Обозначим через $N_б(t)$ – число белых шаров в момент времени $t$, $N_ч(t)$ – число черных шаров в момент времени $t$.

**а)** Покажите, что динамика обратима и $2n$-периодична. Если предположить, что множество $Q$ выбирается случайно, то обратимость теряется. Ограничимся также рассмотрением $t \ll n$. Покажите, что при этих предположениях

$$E_Q[N_б(t) - N_ч(t)] = (1 - 2\mu)^t E_Q[N_б(0) - N_ч(0)].$$

**б)\*** Введем в динамику п. а) случайность. Для этого определим последовательность i.i.d. ($p < 1/2$), независимую от случайности при выборе множества $Q$,

$$\chi(t) = \begin{cases} 1, & \text{с вероятностью } p \\ -1, & \text{с вероятностью } 1-p \end{cases}.$$

Пронумеруем последовательность точек на окружности и сопоставим белому цвету число $+1$, черному $-1$. Введем также функцию (случайную, поскольку $Q$ случайно)

$$\delta_k = \begin{cases} -1, & k \in Q \\ +1, & k \notin Q \end{cases}.$$

Пусть $C_k(t)$ обозначает цвет шара в точке $k$ в момент времени $t$. В таких обозначениях новая динамика будет иметь вид

$$C_k(t+1) = \chi(t+1)\delta_k C_k(t),$$

т.е. с вероятностью $p$ шар еще может (случайно) изменить цвет (все эти обозначения введены для удобства доказательства). Для простоты, будем считать, что в начальный момент все шары белые. Покажите что при $t \ll n$

$$E_{\delta,\chi}\left[\frac{N_\delta(t) - N_ч(t)}{n}\right] \sim (1-2\mu)^t (1-2p)^t,$$

$$D_{\delta,\chi}\left[\frac{N_\delta(t) - N_ч(t)}{n}\right] \sim \frac{1}{n}(1-2\mu)^{2t}(1-2p)^{2t}.$$

Покажите также, что если множество $Q$ зафиксировано (не случайно), то сомножитель со степенью $(1-2\mu)$ исчезает вместе с ограничением $t \ll n$.

**Замечание.** Более подробно об этой модели написано в старых, но до сих пор актуальных и популярных книгах М. Каца, переведенных на русский язык, а также в монографии *Опойцев В.И.* Нелинейная системостатика. М.: Наука, 1986. Отметим, что если перейти при задании множества $Q$ по принципу случайного выбора $m$ позиций, к случайному разыгрыванию каждой позиции: с вероятностью $\mu$ каждая точка (независимо от остальных) принадлежит $Q$, а с вероятностью $1-\mu$ не принадлежит, то получение фактора $(1-2\mu)^t$ существенно упрощается. В статистической физике такой переход и его различные обобщения называют методом большого канонического ансамбля. Хотя сама по себе модель Каца в каком-то смысле является "карикатурной" моделью статистической физики,[20] она, по-прежнему, вызывает интерес ведущих специалистов в этой области и смежных областях (см., например, работы В.В. Козлова, Аджиева–Веденяпина).

---

[20] Именно так, по-видимому, ее и задумывал М. Кац, пытаясь немного обобщить модель из парадокса Эренфестов, см. задачу 1). Тем не менее, эта модель оказалась очень важной стартовой площадкой для многих исследований (в виду своей простоты с одной стороны и возможностью продемонстрировать разнообразные сложности статистической физики с другой: парадокс обратимости, парадокс времени возвращения).